# ON THE CONVERGENCE ANALYSIS OF ASYNCHRONOUS DISTRIBUTED QUADRATIC PROGRAMMING VIA DUAL DECOMPOSITION

KOOKTAE LEE* AND RAKTIM BHATTACHARYA*

**Abstract.** In this paper, we analyze the convergence as well as the rate of convergence of asynchronous distributed quadratic programming (QP) with dual decomposition technique. In general, distributed optimization requires synchronization of data at each iteration step due to the interdependency of data. This synchronization latency may incur a large amount of waiting time caused by an idle process during computation. We aim to attack this synchronization penalty in distributed QP problems by implementing asynchronous update of dual variable. The price to pay for adopting asynchronous computing algorithms is unpredictability of the solution, resulting in a tradeoff between speedup and accuracy. Thus, the convergence to an optimal solution is not guaranteed owing to the stochastic behavior of asynchrony. In this paper, we employ the switched system framework as an analysis tool to investigate the convergence of asynchronous distributed QP. This switched system will facilitate analysis on asynchronous distributed QP with dual decomposition, providing necessary and sufficient conditions for the mean square convergence. Also, we provide an analytic expression for the rate of convergence through the switched system, which enables performance analysis of asynchronous algorithms as compared with synchronous case. To verify the validity of the proposed methods, numerical examples are presented with an implementation of asynchronous parallel QP using OpenMP.

## 1. Introduction.

Recent advancement of distributed and parallel computing technologies has brought massive processing capabilities in solving large-scale optimization problems. Distributed and parallel computing may reduce computation time to find an optimal solution by leveraging the parallel processing in computation. Particularly, distributed optimization will likely be considered as a key element for large-scale statistics and machine learning problems, currently represented by the word "big data". One of the reasons for the preference of distributed optimization in big data is that the size of data set is so huge that each data set is desirably stored in a distributed manner. Thus, global objective is achieved in conjunction with local objective functions assigned to each distributed node, which requires communication between distributed nodes in order to attain an optimal solution.

For several decades, there have been remarkable studies that have enabled to find an optimal solution in a decentralized fashion, for example, dual decomposition [9], [2], [12], [19], [3], augmented Largrangian methods for constrained optimization [21], [29], [16], [1], alternating direction method of multipliers (ADMM) [20], [18], [17], Spingarn's method, [30], Bregman iterative algorithms for $\ell_1$ problems [6], [8], [11], Douglas-Rachford splitting [10], [27], and proximal methods [31]. More details about history of developments on the methods listed above can be found in the literature [5]. In this study, we mainly focus on the analysis of *asynchronous* distributed optimization problems. In particular, we aim to investigate the behavior of asynchrony in the Lagrangian dual decomposition method for *distributed quadratic programming (QP) problems*, where QP problems refer to the optimization problems with a quadratic objective function associated with linear constraints. This type of

[1]Department of Aerospace Engineering, Texas A&M University, College Station, TX 77843-3141 USA (animodor@tamu.edu, raktim@tamu.edu).

QP problems has broad applications including least square with linear constraints, regression analysis and statistics, SVMs, lasso, portfolio optimization problems, etc. With an implementation of Lagrangian dual decomposition, the original QP problems that are separable can be solved in a distributed sense. For this dual decomposition technique, we will study how the asynchronous computing algorithms affect on the convergence as well as the rate of convergence for the dual variable.

Typically, distributed optimization requires synchronization of the data set at each iteration step due to the interdependency of data. For massive parallelism, this synchronization may result in a large amount of waiting time as load imbalance between distributed computing resources would take place at each iteration step. In this case, some nodes that have completed their tasks should wait for others to finish assigned jobs, which causes idle process of computing resources, incurring waste of computation time. In this paper, we attack this restriction on synchronization penalty necessarily required in distributed and parallel computing, through the implementation of *asynchronous computing algorithms*. The asynchronous computing algorithms that do not suffer from synchronization latency thus have a potential to break through the paradigm of distributed and parallel optimization. Unfortunately, it is not completely revealed yet what is the effect of asynchrony on the convergence as well as the rate of that in the distributed optimization. Due to the stochastic behavior of asynchrony, the solution for the asynchronous distributed QP may diverge even if it is guaranteed that the synchronous scheme provides a convergence to an optimal solution. Although Bertsekas [4] introduced a sufficient condition for the convergence of general asynchronous fixed-point iterations (see chapter 6.2), which is equivalent to a diagonal dominance condition for QP problems, however, this condition is known to be very strong and thus conservative, according to the literature [28]. Therefore, the primal emphasis of this research is placed on: 1) convergence analysis; 2) analytic estimation on the rate of convergence, by employing a new framework for analysis of distributed QP problems with an asynchronous update of dual variable.

For this purpose, we will adopt the *switched system* [15], [14], [13], [25], [22], [26], [24], [23] framework as an analysis tool. In general, the switched system is defined as a dynamical system that consists of a set of subsystem dynamics and a certain switching logic that governs a switching between subsystems. For asynchronous algorithms of which dynamics is modeled by the switched system, subsystem dynamics denotes all possible asynchronous computing due to the difference of data processing time in each distributed computing devices. Then, a certain switching logic can be implemented to stand for a random switching between subsystem dynamics. Thus, the switched system framework can be used to properly model the dynamics of asynchronous computing algorithms. Lee *et al.* [24], for example, introduced the switched system to represent the behavior of asynchrony in massively parallel numerical algorithms. In this literature, the authors applied the switched dynamical system framework in order to analyze the convergence, rate of convergence, and error probability for asynchronous parallel numerical algorithms. Based on this switched system framework, this paper provides a new approach for convergence analysis of asynchronous distributed QP problems with dual decomposition technique. The proposed methods will guarantee the convergence to the optimal solution in the mean square sense. In addition, we will study how fast each scheme (e.g., synchronous and asynchronous scheme) converges to an optimal solution by studying the rate of convergence in analytic form. Therefore, this paper will present fundamental yet important analysis on the asynchronous distributed QP problems through the switched system framework, which facilitates investigation on the stochastic behavior of asynchrony.

Rest of this paper is organized as follows. In section 2, preliminaries are presented in connection with problem formulations for asynchronous distributed QP problems using dual decomposition. Section 3 introduces the switched system to model the asynchrony in the asynchronous distributed QP problems. The results for the convergence and the rate of convergence by employing the switched system framework are derived in section 4 and 5, respectively. The numerical example with a real implementation of distributed and parallel QP is provided in section 6, to verify the validity of the proposed methods. Finally, section 7 concludes the paper.

## 2. Preliminaries and Problem Formulation.

**Notation:** The real number, positive integer, and the non-negative integer are denoted by the symbol $\mathbb{R}$, $\mathbb{N}$, and $\mathbb{N}_0$, respectively. The symbol $^\top$ represents the transpose operator. For any real matrix $A, B \in \mathbb{R}^{n\times n}$, the inequality $A < B$ is interpreted by the quadratic sense. (i.e., $v^\top A v < v^\top B v$ for any real vector $v \in \mathbb{R}^n$). In addition, the symbol $\otimes$ stands for the Kronecker product.

### 2.1. Duality Problem.

Consider the following QP problem with a linear inequality constraint.

$$\begin{aligned} &\text{minimize} && f(x) && (2.1)\\ &\text{subject to} && Ax \le b, && (2.2)\end{aligned}$$

where $f(x)$ is given by a quadratic form, meaning $f(x) = \frac{1}{2}x^\top Q x + c^\top x$, the matrix $Q \in \mathbb{R}^{n\times n}$ is a symmetric, positive definite and $c \in \mathbb{R}^n$ is a vector. Further, in the inequality constraint (2.2), it is such that $A \in \mathbb{R}^{m\times n}$ and $b \in \mathbb{R}^m$. If we define the Lagrangian as $L(x, y) \triangleq f(x) + y^\top(Ax - b)$, where $y \in \mathbb{R}^m$ is the dual variable or Lagrange multiplier, then the dual problem for above QP is formulated as follows.

Duality using Lagrangian:

$$\begin{aligned} &\text{maximize} && \inf_x L(x, y) && (2.3)\\ &\text{subject to} && y \ge 0. && (2.4)\end{aligned}$$

The primal optimal point $x^\star$ is obtained from a dual optimal point $y^\star$ as

$$x^\star = \operatorname*{argmin}_x L(x, y^\star).$$

By implementing gradient ascent, one can solve the dual problem, provided that $\inf L(x, y)$ is differentiable. In this case, the iteration to find the $x^\star$ is constructed as follows:

$$x^{k+1} := \operatorname*{argmin}_x L(x, y^k), \tag{2.5}$$

$$y^{k+1} := y^k + \alpha^k(Ax^{k+1} - b), \tag{2.6}$$

where $\alpha^k$ is a step size and the upper script denotes the discrete-time index for iteration.

For the quadratic objective function $f(x)$, the value $\operatorname*{argmin}_x L(x, y^k)$ can be alternatively obtained by $\nabla_x L(x, y^k) = 0$, which leads to

$$\begin{aligned}\operatorname*{argmin}_x L(x, y^k) &= \nabla_x \left( \frac{1}{2} x^\top Q x + c^\top x + {y^k}^\top (Ax - b) \right) \\ &= Qx + c + A^\top y^k = 0.\end{aligned}$$

From (2.5), we have

$$x^{k+1} = -Q^{-1}(A^\top y^k + c). \tag{2.7}$$

Plugging (2.7) into (2.6) results in

$$\begin{aligned}y^{k+1} &= y^k + \alpha^k \left( A \left( -Q^{-1}(A^\top y^k + c) \right) - b \right) \\ &= (I - \alpha^k A Q^{-1} A^\top) y^k - \alpha^k (A Q^{-1} c + b). \end{aligned} \tag{2.8}$$

With the assumption that $y^k \geq 0 \ \forall k$, the above equation provides the solution for $y^\star$ and hence $x^\star$, if $\rho(I - \alpha^k A Q^{-1} A^\top) < 1$ as follows:

$$y^\star = (I - \alpha^k A Q^{-1} A^\top) y^\star - \alpha^k (A Q^{-1} c + b).$$

$$\Rightarrow \begin{cases} y^\star = -(A Q^{-1} A^\top)^{-1} (A Q^{-1} c + b), & \text{(if } A Q^{-1} A^\top \text{ is non-singular)}, \\ x^\star = -Q^{-1}(A^\top y^\star + c). \end{cases} \tag{2.9}$$

**2.2. Dual Decomposition with Synchronous update.** In this subsection, we consider that $f(x) = \frac{1}{2} x^\top Q x + c^\top x$ is *separable*, which means

$$\begin{aligned}f(x) &= \sum_{i=1}^{N} f_i(x_i) \\ &= \sum_{i=1}^{N} \left( \frac{1}{2} x_i^\top Q_i x_i + c_i^\top x_i \right),\end{aligned}$$

where $x = [x_1^\top, x_2^\top, \ldots, x_N^\top]^\top$ and the variables $x_i \in \mathbb{R}^{n_i}$, $i = 1, 2, \ldots, N$ are subvectors of $x$. Also, the matrix $A$ in (2.2) satisfies $Ax = \sum_{i=1}^{N} A_i x_i$, where $A_i$ is such that $A = [A_1, A_2, \ldots, A_N]$.

Then, the equations (2.5) and (2.6) are updated by

$$x_i^{k+1} := \operatorname*{argmin}_{x_i} L(x_i, y^k) = -Q_i^{-1}(A_i^\top y^k + c), \tag{2.10}$$

$$y^{k+1} := y^k + \alpha^k (A x^{k+1} - b). \tag{2.11}$$

Note that when updating $x_i^{k+1}$, $i = 1, 2, \ldots, N$, each value is computed by distributed nodes. Hence, the computation for $x_i^{k+1}$ can be processed in parallel and then, each value of $x_i^{k+1}$ is transmitted to the master node to compute $y^{k+1}$ in the gathering stage. Therefore, as in (2.11), updating $y^{k+1}$ requires synchronization of $x_i^{k+1}$ across all spatial index $i$ at time $k+1$ because $x^{k+1}$ is obtained by stacking $x_i^{k+1}$ from $i = 1$ to $N$. In Fig. 1, we described the conceptual schematic of synchronous update for dual variable $y$. If computing delay occurs among one of the index $i$ due to the difference of processing time in distributed node, the process to update $y^{k+1}$ has to be paused until all data is received from distributed nodes. This implies that the more parallel computing we have, the more delays may take place, resulting in a large amount of the idle time. Consequently, this idle time for synchronization becomes dominant compared to the pure computation time to solve the QP problem in parallel. In massive parallel computing algorithm, it has been reported that the synchronization latency may be up to 50% of total computation time according to the literature [7]. In order to mitigate or avoid this type of restriction that severely affects on the performance to obtain an optimal solution, we introduce *asynchronous computing algorithm* in the following subsection.

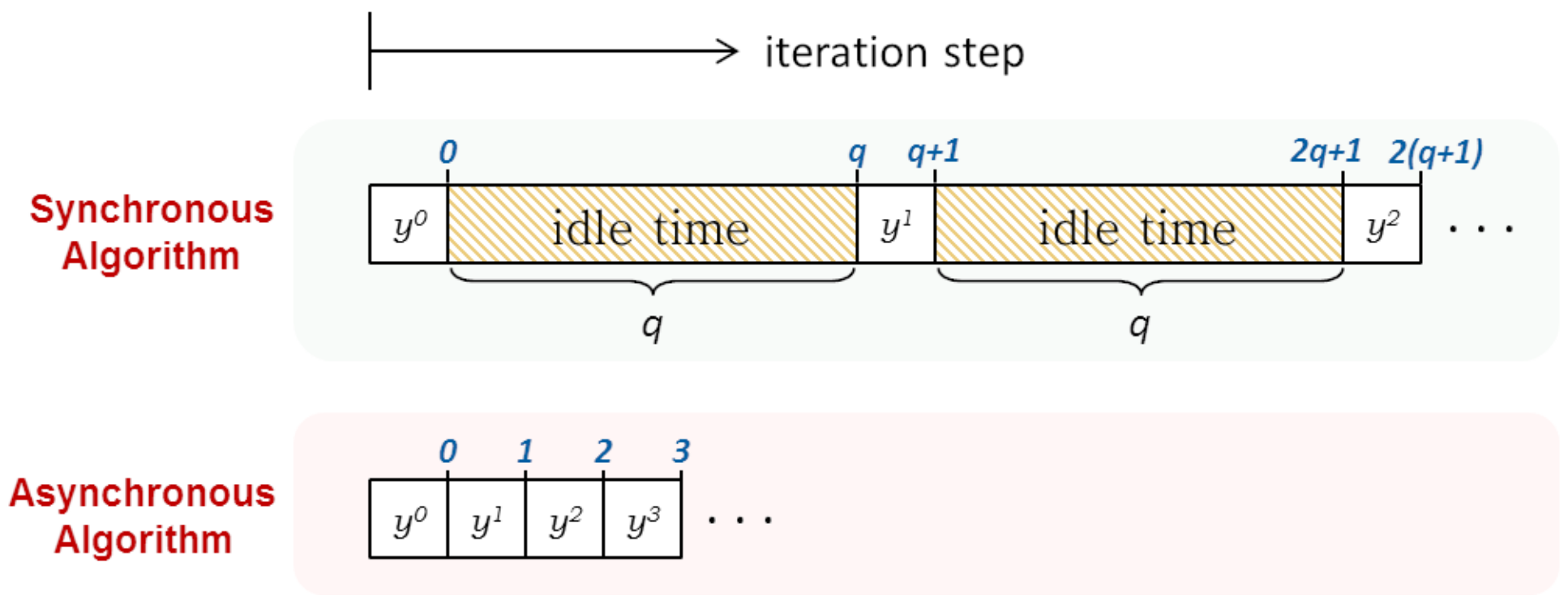

FIG. 1. *The schematic of update timing for the variable $y^k$; upper one shows the synchronous algorithm, where $q$ is the length of maximum delay – i.e., all delays are bounded by $q$; bottom one shows asynchronous algorithm. The time to compute $y^k$ is given by 1 CPU time.*

**2.3. Dual Decomposition with Asynchronous update.** In order to alleviate this synchronization penalty, we consider asynchronous update of dual variable $y$. In this case, the master node to compute $y^{k+1}$ does not wait until all $x_i^{k+1}$ is gathered. Rather, it proceeds with the value for $x_i$ saved in the buffer memory. Thus, $y$ value is updated asynchronously. To model the asynchronous dynamics of dual decomposition, we consider the new state vectors as follows.

- The state for the *Asynchronous* model:

$$\tilde{x}^k := [x_1^{k_1^*\top}, x_2^{k_2^*\top}, \dots, x_N^{k_N^*\top}]^\top,$$

  where $k_i^* \in \{k, k-1, \dots, k-q+1\}$, $i = 1, 2, \dots, N$, denotes delay term that may take place due to the load imbalance in distributed nodes, and the term $q \in \mathbb{N}$ represents the maximum possible delay.

For this asynchronous case, $y$-update is given by

$$\begin{aligned} y^{k+1} &:= y^k + \alpha^k (A\tilde{x}^{k+1} - b) \\ &= y^k + \sum_{i=1}^{N} \left( \alpha_i^k A_i \tilde{x}_i^{k+1} - \frac{1}{N} \alpha_i^k b \right), \end{aligned} \tag{2.12}$$

where $\alpha_i^k$ is the step size for the index $i$.

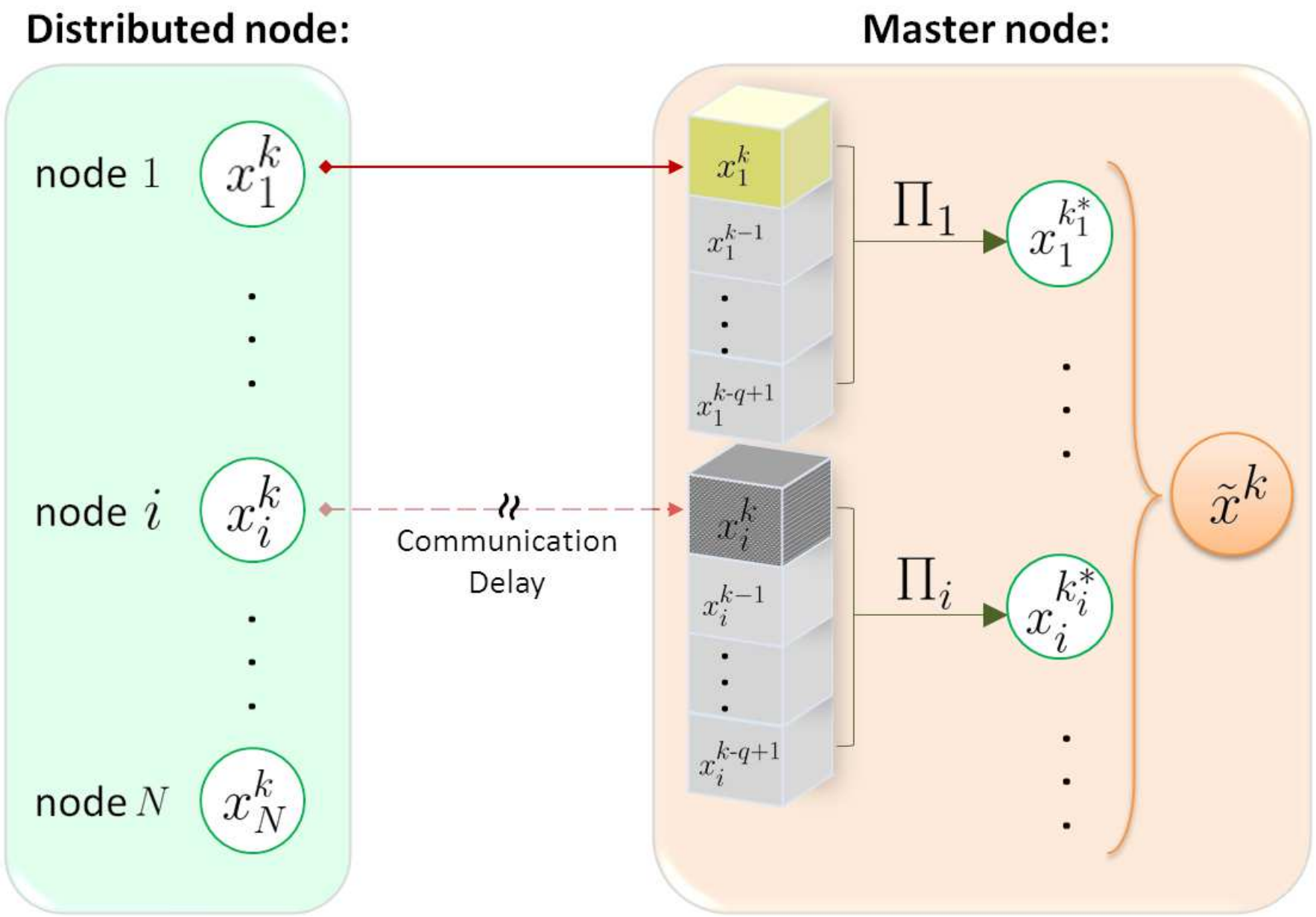

$q$ : the maximum possible delay

$x_i^k$ : the value of $x_i$ at time $k$

$x_i^{k_i^*}$ : the random variable such that $x_i^{k_i^*} \in \{x^k, x^{k-1}, \ldots, x^{k-q+1}\}$

$\Pi_i := [(\pi_1)_i, (\pi_2)_i, \ldots, (\pi_q)_i]$, where $(\pi_j)_i$, $j = 1, 2, \ldots, q$, stands for the modal probability for $x_i^{k_i^*}$

$\tilde{x}^k := [x_1^{k_1^*\top}, x_2^{k_2^*\top}, \ldots, x_N^{k_N^*\top}]^\top$

Fig. 2. *The schematic of the stochatic asynchronous algorithm in the distributed quadratic programming. In this figure, the maximum delay is bounded by $k - q + 1 \leq k_i^* \leq k$, $\forall i$. Each node has the probability $\Pi_i$ to represent random delays.*

Although $\alpha_i^k$ may vary at each time step, we let $\alpha_i^k$ be a constant value, denoted by $\alpha_i$, for simplicity. Hence, it satisfies that $\alpha := \sum_{i=1}^N \alpha_i$, which is a fixed value.

There are two different ways to update dual variable $y$. Throughout the paper, we denote these two different cases as the *deterministic asynchronous algorithm* and the *stochastic asynchronous algorithm*, respectively, in order to clarify and differentiate them. The deterministic asynchronous algorithm stands for the case where the variable $k_i^*$ is considered as a constant value and is given by $k_i^* := k - q + 1$, $\forall i$. Thus, it leads to $\tilde{x}^k := [x_1^{k-q+1\top}, x_2^{k-q+1\top}, \ldots, x_N^{k-q+1\top}]$. In this case, it is assumed that the value $x_i^{k-q+1}$, which is a $q$-step prior value of $x_i^k$, is always available to the master node. In other words, all delays are assumed to be bounded by the finite value $q$. Therefore, one can proceed with $y$-update, given in (2.12), without synchronization when applying the deterministic asynchronous algorithm. Note that there is no randomness in the deterministic asynchronous algorithm. Although this deterministic case obviates the unnecessary idle time by avoiding synchronization, it always utilizes $q$-step prior values saved in the buffer memory. In the real implementation of the distributed optimization, however, $k_i^*$ varies from distributed nodes and also changes over each iteration. Thus, we consider another case by letting $x_i^{k_i^*}$ as a random vector, where $k_i^*$ becomes one of the values in the given set $\{k, k-1, \ldots, k-q+1\}$. To distinguish this case with the deterministic asynchronous algorithm, it is referred to as the stochastic asynchronous algorithm.

Fig. 2 describes the conceptual schematic of the stochastic asynchronous algorithm using the dual decomposition in QP problem. Depending on the processing capability and load balance in distributed nodes, the value for $x_i^k$ is available or not in the master node at each iteration step. We assume that this delay is bounded by the finite value $q$. To describe the randomness of such delays, we adopt a probability $\Pi_i := [(\pi_1)_i, (\pi_2)_i, \ldots, (\pi_q)_i] \in \mathbb{R}^{1\times q}$ that predicts which value for $x_i^k$ will be used to update $y^k$ as shown in Fig. 2.

Starting from (2.12), with the definition of the set $\mathcal{S}^k := \{k_i^* | k_i^* = k\}$ and the symbol $\Phi_i := -\alpha_i A_i Q_i^{-1} A_i^\top$, the state dynamics of the stochastic asynchronous algorithm is then given by

$$\begin{aligned}
y^{k+1} &= y^k + \sum_{i=1}^{N} \left( \alpha_i A_i \tilde{x}_i^{k+1} - \frac{1}{N} \alpha_i^k b \right) \\
&= y^k + \sum_{i \in \mathcal{S}^{k+1}} \alpha_i A_i x^{k+1} + \sum_{i \in \mathcal{S}^k} \alpha_i A_i x^k + \cdots \\
&\qquad\qquad + \sum_{i \in \mathcal{S}^{k-q+2}} \alpha_i A_i Q_i^{-1} A_i^\top x^{k-q+2} + \left( \sum_{i=1}^{N} -\frac{1}{N} \alpha_i b \right) \\
&= \left( I - \sum_{i \in \mathcal{S}^{k+1}} \Phi_i \right) y^k - \left( \sum_{i \in \mathcal{S}^k} \Phi_i \right) y^{k-1} - \cdots \qquad \text{(by (2.10))} \\
&\qquad - \left( \sum_{i \in \mathcal{S}^{k-q+2}} \Phi_i \right) y^{k-q+1} + \left( \sum_{i=1}^{N} -\alpha_i A_i Q_i^{-1} c - \frac{1}{N} \alpha_i b \right). \qquad (2.13)
\end{aligned}$$

The above equation is simplified by the following definitions, given by

$$R_i(k) := \sum_{j \in \mathcal{S}^{k-i+2}} \Phi_j, \qquad (2.14)$$

$$B := \left( \sum_{i=1}^{N} -\alpha_i A_i Q_i^{-1} c - \frac{1}{N} \alpha_i b \right), \qquad (2.15)$$

resulting in

$$y^{k+1} = (I - R_1(k))\, y^k - R_2(k) y^{k-1} - \cdots - R_q(k) y^{k-q+1} + B, \qquad (2.16)$$

where the time-varying matrix $R_i(k)$ completely depends on the value $k_i^*$ that is a random event.

As described in [4], it is a very challenging task to analyze the stochastic asynchronous algorithm (see page 101, chapter 1). The primary goal of this paper is, therefore, to analyze not only the convergence but also the rate of that for the stochastic asynchronous algorithm which brings stochastic process for the state $y^k$. For this purpose, we adopt a *switched linear system (or jump linear system, interchangeably)* framework that will be introduced in the next section in more detail.

**3. A Switched System Approach for Asynchronous Computing Algorithms.** In order to solve the dual decomposition problem with random delays in distributed nodes, we define a new augmented state $Y^k := [{y^k}^\top, {y^{k-1}}^\top, \ldots, {y^{k-q+1}}^\top]^\top$. Then, one can define the following recursive dynamics:

$$
\underbrace{\begin{bmatrix} y^{k+1} \\ y^k \\ y^{k-1} \\ \vdots \\ y^{k-q+2} \end{bmatrix}}_{=Y^{k+1}} = \underbrace{\begin{bmatrix} I - R_1(k) & -R_2(k) & -R_3(k) & \cdots & -R_q(k) \\ I & 0 & \cdots & & 0 \\ 0 & I & 0 & \cdots & 0 \\ \vdots & & & \ddots & \vdots \\ 0 & 0 & & I & 0 \end{bmatrix}}_{=W(k)} \underbrace{\begin{bmatrix} y^k \\ y^{k-1} \\ y^{k-2} \\ \vdots \\ y^{k-q+1} \end{bmatrix}}_{=Y^k} + \underbrace{\begin{bmatrix} B \\ 0 \\ 0 \\ \vdots \\ 0 \end{bmatrix}}_{=C}, \tag{3.1}
$$

where $I$ and $0$ are identity and zero matrices with proper dimensions, respectively. Consequently, the above recursive equation ends up with the following simple form:

$$
\Rightarrow Y^{k+1} = W(k)Y^k + C
$$

In fact, the structure of the time-varying matrix $W(k)$ is not arbitrary, but it has a finite number of forms, given by $q^N$, which counts all possible scenarios to distribute $N$ numbers of $\Phi_i$, $i = 1, 2, \ldots, N$, matrices into the finite number of $q$. In the switched system, this number is referred to as the *"switching mode number"*, and we particularly denote this number with the symbol $m$. For instance, when $q = 2$ and $N = 2$, the switching mode number is given by $m = 2^2 = 4$. Thus, at each time $k$, the matrix $W(k)$ has one of the following form:

$$
W_1 = \begin{bmatrix} I - \Phi_1 - \Phi_2 & 0 \\ I & 0 \end{bmatrix}, \qquad W_2 = \begin{bmatrix} I - \Phi_1 & -\Phi_2 \\ I & 0 \end{bmatrix},
$$
$$
W_3 = \begin{bmatrix} I - \Phi_2 & -\Phi_1 \\ I & 0 \end{bmatrix}, \qquad W_4 = \begin{bmatrix} I & -\Phi_1 - \Phi_2 \\ I & 0 \end{bmatrix}.
$$

Then, only one out of all set of matrices $\{W_r\}_{r=1}^m$ will be used at each time $k$ to update the system state $Y^k$, which results in the *switched linear system* structure as follows.

Consider the switched system:

$$
Y^{k+1} = W_{\sigma_k} Y^k + C, \quad \sigma_k \in \{1, 2, \ldots, m\},\ k \in \mathbb{N}_0, \tag{3.2}
$$

where $\{\sigma_k\}_{k=0}^\infty$ denotes the switching sequence that describes how the asynchrony takes place. Then, the switching probability $\Pi(k) := \Pi_1(k) \otimes \Pi_2(k) \otimes \cdots \otimes \Pi_N(k) = [\pi_1(k), \pi_2(k), \ldots, \pi_m(k)]$, where $\Pi_i(k)$ represents the probability for $x_i^{k_i^*}$ as depicted by Fig. 2, determines which mode $\sigma_k$ will be utilized at each time step. (Note that $\Pi_i(k)$ and hence $\Pi(k)$ are not necessarily to be stationary.) In this case, the switched linear system is named by "stochastic switched linear system" or "stochastic jump linear system" [26] because the switching is a stochastic process. The benefit when applying this stochastic switched linear system structure is that the delay in the asynchronous algorithm is naturally taken into account by the switched system framework. Hence, the randomness of the asynchronous algorithm is represented by a certain switching logic.

REMARK 3.1. ***(Computational complexity due to an extremely large number of the switching modes)*** *Although the stochastic switched linear system framework is suitable for modeling the dynamics of the stochastic asynchronous algorithm in distributed QP problems, it results in an extremely large number of the switching modes, causing computational complexity. For instance, even if $q = 2$ and $N = 20$, we have $m = q^N = 2^{20}$, and it is impractical to store such large numbers of matrices in the real implementation. Therefore, it is necessary to develop proper methods to analyze the stochastic asynchronous algorithm using the switched linear system without any concerns for such computational complexity issues.*

To avoid the computational complexity problems stated above, we firstly make following assumptions for analysis of both the convergence and the rate of the convergence for the stochastic asynchronous algorithm:

- **Assumption 3.1.** We consider the random delays that occur during the computation of $x_i^k$ at each node. In this case, the probability $\Pi_i(k) = [(\pi_1(k))_i, (\pi_2(k))_i, \ldots, (\pi_q(k))_i]$ describes which value for $x_i^{k_i^*}$ will be used among the given set $\{x_i^k, x_i^{k-1}, \ldots, x_i^{k-q+1}\}$. Then, we assume that each modal probability $(\pi_j(k))_i$ is *stationary*, and hence $\Pi_i(k)$ is also stationary in time.

Under the Assumption 3.1., the switching probability $\Pi(k) := \Pi_1 \otimes \Pi_2 \otimes \cdots \otimes \Pi_N$ becomes stationary. For this case, the jump linear system with the given dynamics in (3.2) is termed as the independent, identically distributed (i.i.d.) jump linear system. Since the modal switching probability $\pi_r$ is a probability, it satisfies $0 \leq \pi_r \leq 1$, $\forall r$ and $\sum_{r=1}^{m} \pi_r = 1$. This stationary occupation probability rules which system matrix $W_r$ will be used at each instance. The implementation of the switching sequence $\{\sigma_k\}$, governed by $\Pi$, describes the randomness for the stochastic asynchronous algorithm *in an average sense.*

**4. Convergence Analysis.** In this section, the convergence of the state $Y^k$ for the stochastic asynchronous model will be studied under the switched system framework. For several decades, the stability results for the switched systems with stochastic jumping parameters have been well established, for example, in the literature [25], [26], [15]. However, these methods are inapplicable to the asynchronous computing algorithm with massive parallelism because it results in extremely large numbers of switching modes, leading to *computational complexity* as explained in Remark 3.1. Therefore, we aim to investigate the convergence and the rate of convergence for the asynchronous algorithm without any concerns for such computational complexity issues. Particularly, this section will provide a convergence condition for the stochastic asynchronous algorithm in distributed QP problems.

Before proceeding further to investigate the asynchronous model, we analyze the convergence of the *synchronous* case *without delays* for a reference. Since in the synchronous algorithm all values are synchronized after each iteration, no delays occur when updating the state dynamics. Then, the state $Y_{\text{sync.}}^k$ for the synchronous case is governed by the following recursive equation:

$$
(4.1)\qquad \underbrace{\begin{bmatrix} y^{k+1} \\ y^{k} \\ y^{k-1} \\ \vdots \\ y^{k-q+2} \end{bmatrix}}_{=Y_{\text{sync.}}^{k+1}} = \underbrace{\begin{bmatrix} I-R & 0 & 0 & \cdots & 0 \\ I & 0 & \cdots & & 0 \\ 0 & I & 0 & \cdots & 0 \\ \vdots & & \ddots & & \vdots \\ 0 & 0 & & I & 0 \end{bmatrix}}_{=W_{\text{sync.}}} \underbrace{\begin{bmatrix} y^{k} \\ y^{k-1} \\ y^{k-2} \\ \vdots \\ y^{k-q+1} \end{bmatrix}}_{=Y_{\text{sync.}}^{k}} + \underbrace{\begin{bmatrix} B \\ 0 \\ 0 \\ \vdots \\ 0 \end{bmatrix}}_{=C},
$$

where the matrix $R := \sum_{i=1}^{q} R_i(k) = \sum_{j=1}^{N} \Phi_j$ is time-invariant, and hence the matrix $W_{\text{sync.}}$ is also constant. Then, the steady-state value of $Y_{\text{sync.}}^{\star} := \lim_{k\to\infty} Y_{\text{sync.}}^{k}$, is obtained by

$$
(4.2)\qquad \begin{aligned} & Y_{\text{sync.}}^{\star} = W_{\text{sync.}} Y_{\text{sync.}}^{\star} + C. \\ \Rightarrow\ & Y_{\text{sync.}}^{\star} = (I - W_{\text{sync.}})^{-1} C, \qquad \big(\text{if } (I - W_{\text{sync.}}) \text{ is non-singular}\big) \end{aligned}
$$

if the condition $\rho(W_{\text{sync.}}) < 1$ holds.

However, the state in the i.i.d. switched linear system that represents the stochastic asynchronous model, evolves with the dynamics given in (3.2), where the matrix $W_{\sigma_k}$ is determined by the switching probability $\Pi$. Thus, the state of the asynchronous model becomes a random vector, obstructing the convergence analysis of the stochastic asynchronous model. For the stochastic switched systems, various convergence (stability) notions have been developed [15], to guarantee the system stability. Among different convergence notions, we will focus on the mean square convergence, defined below.

Definition 4.1. *(Definition 1.1, [13]) The switched system is said to be mean square stable (convergent) if for any initial condition $x_0$ and arbitrary initial probability distribution $\Pi(0)$, $\lim_{k\to\infty} \mathbb{E}\big[||x(k,x_0) - x^\star||^2\big] = 0$, where $x^\star$ is the fixed-point value of $x^k$, i.e. $\lim\limits_{k\to\infty} x^k = x^\star$.*

The necessary and sufficient condition for the mean square convergence of the i.i.d. jump linear systems is described as follows:

Proposition 4.2. *(Corollary 2.7, [14]) Consider an i.i.d. jump linear system, where $\Pi(k)$ is a stationary probability vector $\{\pi_1, \pi_2, \cdots, \pi_m\}$ for all $k$. Then, the i.i.d. jump linear system is mean square stable (convergent) if and only if the matrix $\sum_{j=1}^{m} \pi_j (W_j \otimes W_j)$ is Schur stable, i.e.*

$$
(4.3)\qquad \rho\left(\sum_{j=1}^{m} \pi_j \left(W_j \otimes W_j\right)\right) < 1.
$$

Once again, massive parallelism results in large $m$, causing computational intractability. Thus, implementation of Proposition 4.2 is unfeasible to analysis of asynchronous distributed and parallel QP problems with massively parallel computing algorithm because the equation in (4.3) requires the summation over index $i$ from 1 up to $m$. In order to avoid this problem, we provide Algorithm 1.

By executing Algorithm 1 at every time step in the master node, the random vector $\tilde{x}^k$ has the following form: $\tilde{x}^k = [(x_1^{\xi})^\top, (x_2^{\xi})^\top, \ldots, (x_N^{\xi})^\top]^\top$, where $\xi$ denotes the oldest time among the recently updated values across the index $i = 1, 2, \ldots, N$.

**Algorithm 1**

1: $k_i^* \leftarrow$ one of the values in $\{k, k-1, \ldots, k-q+1\}$ with probability $\Pi_i$.
2: $\xi \leftarrow k$
3: **for** $i \leq N$ **do**
4:     **if** $\xi \leq k_i^*$ **then**
5:         $\xi \leftarrow k_i^*$.
6:         $i \leftarrow i+1$.
7:     **end if**
8: **end for**
9: $\tilde{x}^k \leftarrow [(x_1^{\xi})^\top, (x_2^{\xi})^\top, \ldots, (x_N^{\xi})^\top]^\top$

For example, if $k_i^* = k-2$ for some $i$ is the oldest value over all $k_i^*$, $i = 1, 2, \ldots, N$, then we have $\tilde{x}^k = [(x_1^{k-2})^\top, (x_2^{k-2})^\top, \ldots, (x_N^{k-2})^\top]^\top$. In this case, the modal matrix $W_r$ has the same structure with $W(k)$, given in (3.1), where $R_i(k)$ satisfies

$$R_i(k) = \begin{cases} R, & (\text{if } i = k - \xi + 1) \\ 0. & (\text{otherwise}) \end{cases}$$

The utilization of Algorithms 1 then drastically reduces the switching mode number by $q$ regardless of the value $N$, due to the fact that at each iteration step we intentionally use the oldest updated value saved in buffer memory. For example, when $q = 2$, the matrix $W_{\sigma_k}$ becomes one of the following form:

$$W_1 = \begin{bmatrix} I - R & 0 \\ I & 0 \end{bmatrix}, \qquad W_2 = \begin{bmatrix} I & -R \\ I & 0 \end{bmatrix}.$$

Since Algorithm 1 works as if it aggregates some subsets of the given switching modes, we need to redefine the switching probability $\Pi$ accordingly. Then, $\Pi$ is obtained by the following Theorem.

THEOREM 4.3. *Consider the i.i.d. switched linear system given in* (3.2) *with the switching probability* $\Pi = \Pi_1 \otimes \Pi_2 \otimes \ldots \otimes \Pi_N \in \mathbb{R}^{1 \times q^N}$. *After the implementation of Algorithm 1, the switching probability is redefined by* $\Pi := [\pi_1, \pi_2, \ldots, \pi_q] \in \mathbb{R}^{1 \times q}$, *of which modal probability* $\pi_i$ *has the following form:*

$$\pi_r := \prod_{i=1}^{N} \left( \sum_{j=1}^{r} (\pi_j)_i \right) - \left( \sum_{j=1}^{r-1} \pi_j \right), \quad r = 1, 2, \ldots, q, \tag{4.4}$$

*where the term* $(\pi_j)_i$ *denotes* $j^{th}$ *modal probability for* $\Pi_i$ *(i.e.,* $\Pi_i = [(\pi_1)_i, (\pi_2)_i, \ldots, (\pi_q)_i]$ *).*

*Proof.* For simplicity of the proof, we assume that $N = 2$. The most general case is then proved similarly by induction. In this case, the master node takes the values for each $x_i^{k_i^*}$ according to the probability $\Pi_i$, $i = 1, 2$, which are given by

$$\Pi_1 = [(\pi_1)_1, (\pi_2)_1, \ldots, (\pi_q)_1],$$
$$\Pi_2 = [(\pi_1)_2, (\pi_2)_2, \ldots, (\pi_q)_2].$$

We let the index $j \in \{k, k-1, \ldots, k-q+1\}$ be the value explained in Algorithm 1.

When $j = 1$, the modal switching probability $\pi_1$ is obtained by

$$\begin{aligned}\pi_1 &= \mathbf{Pr}\Big(k_1^* = k,\ k_2^* = k\Big)\\ &= \mathbf{Pr}\Big(k_1^* = k\Big) \times \mathbf{Pr}\Big(k_2^* = k\Big) \qquad \text{(since } k_1^* \text{ and } k_2^* \text{ are independent)}\\ &= (\pi_1)_1 \times (\pi_1)_2.\end{aligned}$$

Similarly, when $j = 2$, we have

$$\begin{aligned}\pi_2 &= \mathbf{Pr}\Big(k_1^* \in \{k, k-1\},\ k_2^* \in \{k, k-1\}\Big) - \pi_1\\ &= \sum_{j=1}^{2} \mathbf{Pr}\Big(k_1^* = k - j + 1\Big) \times \sum_{j=1}^{2} \mathbf{Pr}\Big(k_2^* = k - j + 1\Big) - \pi_1\\ &= \left(\sum_{j=1}^{2} (\pi_j)_1\right) \times \left(\sum_{j=1}^{2} (\pi_j)_2\right) - \pi_1.\end{aligned}$$

In the first line of above equation, we have to extract $\pi_1$ because it corresponds to the case when $j = 1$.

For any arbitrary value $j$ satisfying $j \in \{k, k-1, \ldots, k-q+1\}$, the switching probability is therefore obtained by induction as follows:

$$\begin{aligned}\pi_r &= \mathbf{Pr}\Big(k_1^* \in \{k, k-1, \ldots, k-r+1\},\ k_2^* \in \{k, k-1, \ldots, k-r+1\}\Big) - \sum_{j=1}^{r-1} \pi_j\\ &= \left(\sum_{j=1}^{r} (\pi_j)_1\right) \times \left(\sum_{j=1}^{r} (\pi_j)_2\right) - \sum_{j=1}^{r-1} \pi_j.\end{aligned}$$

Thus, the most general case with $q, N \in \mathbb{N}$ can be induced as follows:

$$\pi_r = \prod_{i=1}^{N} \left(\sum_{j=1}^{r} (\pi_j)_i\right) - \left(\sum_{j=1}^{r-1} \pi_j\right), \quad r = 1, 2, \ldots, q.$$

□

For comparison, the switching mode number without the proposed algorithm is given by $m = q^N$ of which growth is exponential with respect to $N$, whereas with the proposed Algorithm 1, it is given by $m = q$ that is a constant value *irrespective of $N$*. Thus, by leveraging the proposed algorithm, one can apply the mean square convergence condition given in Proposition 4.2, to test the stability of the stochastic asynchronous algorithm. Note that the implementation of Proposition 4.2 was computationally intractable without Algorithm 1 due to the large numbers in $m$. Consequently, the proposed algorithm enables the convergence analysis of the stochastic asynchronous parallel computing algorithm in QP problems.

Once the condition (4.3) is guaranteed with a given i.i.d. switching probability $\Pi$ by implementing Algorithm 1, the steady-state (fixed-point) value $Y^\star := \lim_{k\to\infty} Y^k$, where $Y^k$ is the state for the stochastic asynchronous algorithm of which dynamics is given in (3.2), can be obtained according to Definition 4.1 and is given by

$$\begin{aligned}& Y^\star = W_{\sigma_k} Y^\star + C.\\ \Rightarrow\ & Y^\star = (I - W_{\sigma_k})^{-1} C.\end{aligned} \tag{4.5}$$

Interestingly, $Y^\star$ becomes a unique vector, regardless of $\sigma_k$ that changes over time, due to the inherent structure in matrices $W_{\sigma_k}$ and $C$, which results in $Y^\star = Y^\star_{\text{sync.}}$, where $Y^\star_{\text{sync.}}$ is defined in (4.2). Therefore, the state for the stochastic asynchronous algorithm, denoted by $Y^k$, converges to the unique, identical fixed-point value $Y^\star$, if the condition (4.3) holds.

**5. Rate of Convergence Analysis.** Since the rate of convergence provides information regarding how fast each scheme converges to the fixed-point value, it works as a guideline that suggests which methods will solve the given QP problem faster than other schemes. Therefore, the comparison for the rate of convergence between different schemes is advantageous in terms of estimating the time to obtain an optimal solution for the QP problem. Although asynchronous algorithms are considered to be more time-efficient for obtaining an optimal solution, it is not analytically proved yet what is the rate of convergence. Therefore, in this section we investigate the rate of convergence for three different algorithms (e.g., synchronous, deterministic asynchronous, and stochastic asynchronous algorithms) in analytic form.

**i) Synchronous algorithm with delays:**

For synchronous scheme, $Y^k$ is updated after a certain amount of time due to the idle time for synchronization. As described in Fig. 1, we assume that all data from distributed nodes arrive at the master node within a bounded time $q$. In this case, idle process time for the synchronization is given by $q$ and $Y^k$ can be updated at every $t(q+1)$ time step, where $t \in \mathbb{N}_0$. Consequently, at each time step, $Y^k$-update is given by

$$\begin{aligned}
&\text{at time } t=1: && Y^{(q+1)} = W_{\text{sync.}} Y^0 + C \\
&\text{at time } t=2: && Y^{2(q+1)} = W_{\text{sync.}} Y^{(q+1)} + C \\
&\text{at time } t=3: && Y^{3(q+1)} = W_{\text{sync.}} Y^{2(q+1)} + C \\
&\quad\vdots && \quad\vdots \\
&\text{at arbitrary time } t+1: && Y^{(t+1)(q+1)} = W_{\text{sync.}} Y^{t(q+1)} + C, \quad t \in \mathbb{N}_0
\end{aligned}$$

Now, we consider the term $||Y^k - Y^\star||_\infty$ in order to investigate the rate of convergence for the synchronous algorithm. Then, from the dynamics for synchronous case, given by $Y^k = W_{\text{sync.}} Y^{k-1} + C$, we have

$$\begin{aligned}
||Y^k - Y^\star||_\infty &= ||W_{\text{sync.}} Y^{k-1} + C - Y^\star||_\infty \\
&= ||W_{\text{sync.}} Y^{k-1} - W_{\text{sync.}} Y^\star||_\infty && \big(\text{by } (4.2)\big) \\
&= ||W_{\text{sync.}} \big(W_{\text{sync.}} Y^{k-2} + C\big) - W_{\text{sync.}} Y^\star||_\infty \\
&= ||(W_{\text{sync.}})^2 Y^{k-2} + W_{\text{sync.}}(C - Y^\star)||_\infty \\
&= ||(W_{\text{sync.}})^2 \left(Y^{k-2} - Y^\star\right)||_\infty && \big(\text{by } (4.2)\big) \\
&\quad\vdots \\
&= ||(W_{\text{sync.}})^k \left(Y^0 - Y^\star\right)||_\infty \\
&\le ||(W_{\text{sync.}})^k||_\infty \cdot ||Y^0 - Y^\star||_\infty,
\end{aligned}$$

where $k = t(q+1)$, $t \in \mathbb{N}_0$. Thus, we have the upper bound of the rate of convergence for the synchronous algorithm as follows:

$$||Y^k - Y^\star||_\infty \le ||(W_{\text{sync.}})^k||_\infty \cdot ||Y^0 - Y^\star||_\infty, \quad k = t(q+1),\ t \in \mathbb{N}_0. \tag{5.1}$$

**ii) Deterministic Asynchronous algorithm:**

As described in section 2.3., the deterministic asynchronous algorithm takes advantage of the $q$ step prior value instead of waiting for all $x_i$ values being gathered in the master node for synchronization. In this case, the system dynamics for the deterministic asynchronous scheme is given by

$$Y^{k+1} = W_{\text{det.async.}} Y^k + C,$$

where the matrix $W_{\text{det.async.}}$ is defined as

$$W_{\text{det.async.}} := \begin{bmatrix} I & 0 & 0 & \cdots & -R \\ I & 0 & 0 & \cdots & 0 \\ 0 & I & 0 & \cdots & 0 \\ \vdots & & \ddots & & \vdots \\ 0 & 0 & & I & 0 \end{bmatrix}$$

because in this case we have $\forall i \in \mathcal{S}^{k-q+2}$ in (2.13) for the deterministic asynchronous algorithm, leading to above system dynamics.

Similarly to the process in obtaining (5.1), the upper bound of the rate of convergence for the deterministic asynchronous algorithm is derived by

$$\begin{aligned} ||Y^k - Y^\star||_\infty &= ||(W_{\text{det.async.}})^k \left(Y^0 - Y^\star\right)||_\infty \\ &\leq ||(W_{\text{det.async.}})^k||_\infty \cdot ||Y^0 - Y^\star||_\infty, \quad k \in \mathbb{N}_0. \end{aligned} \tag{5.2}$$

**iii) Stochastic Asynchronous algorithm:**

Since the state $Y^k$ becomes a random vector in the stochastic asynchronous case, the rate of convergence for $||Y^k - Y^\star||_\infty$ forms a distribution rather than a deterministic value, and is difficult to analyze such a distribution. Thus, we take the expectation for $Y^k$ with respect to the i.i.d. switching probability $\Pi$, and investigate the rate of convergence for $||\mathbb{E}[Y^k] - Y^\star||_\infty$.

Under the assumption that the mean square convergence condition in Proposition 4.2 holds, the fixed-point value for $Y^k$ is deterministically given by $Y^\star$. Therefore, it satisfies $\mathbb{E}[Y^\star] = Y^\star$. Taking the expectation in (4.5) results in $\mathbb{E}[Y^\star] = Y^\star = \mathbb{E}[W_{\sigma_k} Y^\star + C] = \mathbb{E}[W_{\sigma_k}]Y^\star + C = \mathbf{Pr}\left(\sum_{r=1}^{q} \pi_r W_r\right) Y^\star + C$. By defining a new matrix $\Lambda := \sum_{r=1}^{q} \pi_r W_r$, we end up with

$$\mathbb{E}[Y^\star] = Y^\star = \Lambda Y^\star + C. \tag{5.3}$$

Then, the term $||\mathbb{E}[Y^k] - Y^\star||_\infty$ becomes

$$
\begin{aligned}
||\mathbb{E}[Y^k] - Y^\star||_\infty &= ||\mathbb{E}[W_{\sigma_{k-1}} Y^{k-1} + C] - Y^\star||_\infty \\
&= ||\sum_{r=1}^{q} \mathbf{Pr}\left(W_{\sigma_{k-1}} Y^{k-1} + C | \sigma_{k-1} = r\right) \underbrace{\mathbf{Pr}\left(\sigma_{k-1} = r\right)}_{=\pi_r} - Y^\star||_\infty \\
&= ||\sum_{r=1}^{q} \pi_r W_r \mathbf{Pr}\left(Y^{k-1} | \sigma_{k-1} = r\right) + C - Y^\star||_\infty \\
&= ||\sum_{r=1}^{q} \pi_r W_r \mathbf{Pr}\left(Y^{k-1} | \sigma_{k-1} = r\right) - \Lambda Y^\star||_\infty \qquad \text{(by (5.3))} \\
&= ||\underbrace{\left(\sum_{r=1}^{q} \pi_r W_r\right)}_{=\Lambda} \sum_{s=1}^{q} \mathbf{Pr}\left(W_{\sigma_{k-2}} Y^{k-2} + C | \sigma_{k-2} = s\right) \pi_s - \Lambda Y^\star||_\infty \\
&= ||\Lambda \left(\sum_{s=1}^{q} \pi_s W_s \mathbf{Pr}\left(Y^{k-2} | \sigma_{k-2} = s\right) + C\right) - \Lambda Y^\star||_\infty \\
&= ||\Lambda \left(\sum_{s=1}^{q} \pi_s W_s \mathbf{Pr}\left(Y^{k-2} | \sigma_{k-2} = s\right)\right) + \Lambda C - \Lambda\left(\Lambda Y^\star + C\right)||_\infty \\
&\qquad \text{(by (5.3))} \\
&= ||\Lambda \left(\sum_{s=1}^{q} \pi_s W_s \mathbf{Pr}\left(Y^{k-2} | \sigma_{k-2} = s\right)\right) - (\Lambda)^2 Y^\star||_\infty \\
&\quad \vdots \\
&= ||(\Lambda)^{k-1} \left(\sum_{t=1}^{q} \pi_t W_t \mathbf{Pr}\left(Y^0 | \sigma_0 = t\right) + C\right) - (\Lambda)^k Y^\star||_\infty \\
&= ||\left(\Lambda\right)^{k-1} \underbrace{\left(\sum_{t=1}^{q} \pi_t W_t\right)}_{=\Lambda} Y_0 - (\Lambda)^k Y^\star||_\infty \\
&= ||(\Lambda)^k \left(Y^0 - Y^\star\right)||_\infty \\
&\leq ||(\Lambda)^k||_\infty \cdot ||Y^0 - Y^\star||_\infty,
\end{aligned}
$$

where we used the law of total probability in above equations.

Therefore, the rate of convergence for the asynchronous scheme is given by:

$$
||\mathbb{E}[Y^k] - Y^\star||_\infty \leq ||(\Lambda)^k||_\infty \cdot ||Y^0 - Y^\star||_\infty, \tag{5.4}
$$

where with implementation of Algorithm 1 the matrix $\Lambda := \sum_{r=1}^{q} \pi_r W_r$ has the following form:

$$
\Lambda = \begin{bmatrix} I - \pi_1 R & -\pi_2 R & -\pi_3 R & \cdots & -\pi_q R \\ I & 0 & 0 & \cdots & 0 \\ 0 & I & 0 & \cdots & 0 \\ \vdots & & \ddots & & \vdots \\ 0 & 0 & \cdots & I & 0 \end{bmatrix}, \quad R := \sum_{i=1}^{N} R_i(k) = \sum_{j=1}^{N} \Phi_j. \tag{5.5}
$$

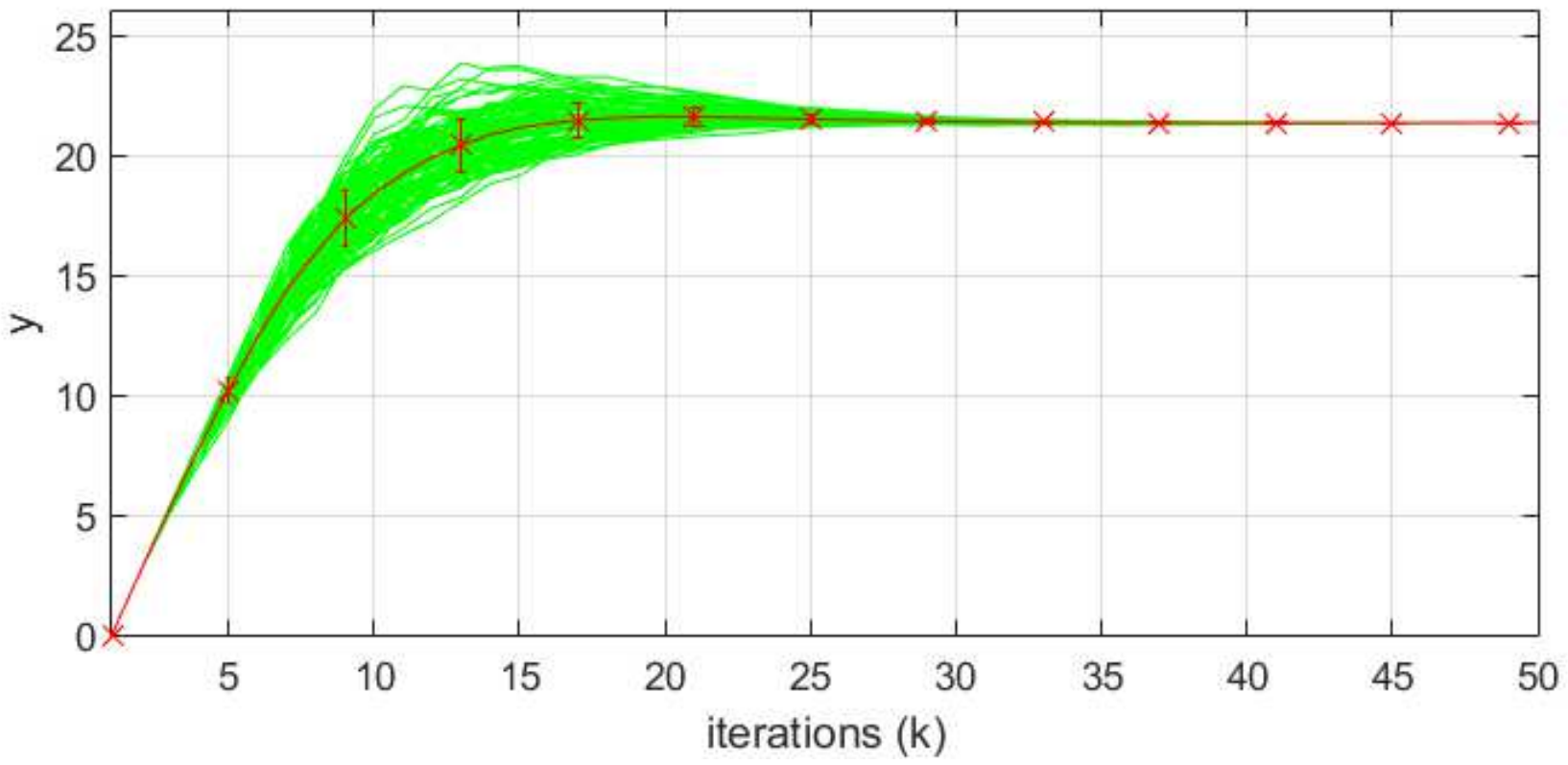

FIG. 3. *Convergence results for distributed quadratic programming with stochastic asynchronous algorithm. The (green) solid lines represent the state trajectory for $y$ with total 100 Monte Carlo simulations (initial value was deterministically given by $y(0) = 2$ for all cases). The (red) solid-cross line denotes the mean and the standard deviation of multiple trajectories, respectively.*

**6. Numerical Example.** In this section, we test the proposed asynchronous algorithms on distributed QP problems with dual decomposition technique. The system for the test bed is given by `Intel(R) Core(TM) i7-4710HQ CPU`, which has 4 cores with 8 threads (by Hyper-Threading Technology), with `8GM memory`. Although, the number of threads for this test bed is not very large, the system is enough to show the performance of proposed asynchronous computing algorithms for distributed QP with dual decomposition. We implemented parallel processing through `OpenMP` API (Application Program Interface) developed for direct multi-threaded, shared memory parallelism.

Let us consider the following distributed QP problem:

$$\begin{array}{ll} \text{minimize} & \displaystyle\sum_{i=1}^{N} \left( \frac{1}{2} x_i^\top Q_i x_i + c_i^\top x_i \right) \\ \text{subject to} & A_i x_i \leq b_i, \qquad\qquad i = 1, 2, \ldots, N. \end{array}$$

The positive definite matrices $Q_i$, the matrices $A_i$, and the vectors $c_i$ and $b_i$ were generated by implementing PSEUDO RANDOM NUMBER generator in `C++`. The dimension of matrices and vectors are set to be: $Q_i \in \mathbb{R}^{n\times n}$, $A_i \in \mathbb{R}^{1\times n}$, $c_i \in \mathbb{R}^{n\times 1}$, and $b_i \in \mathbb{R}$, $i = 1, 2, \ldots, N$, where $n = 10$, $N = 20000$. Thus, computational burden for solving each distributed QP is low, whereas the total number of distributed QP is extremely high. We let the buffer length $q = 8$ and the step size $\alpha_i = 0.27$, $\forall i$.

For this type of massively distributed QP problem, the time for synchronization may become dominant in the total amount of time to solve QP. In this case, asynchronous computing algorithms may lead to speedup by avoiding synchronization. We solved above distributed QP problem with the implementation of the proposed stochastic asynchronous algorithm. In Fig. 3, total 100 times of state trajectories for the dual variable $y$ are given by (green) solid lines. Since $y$-update is stochastic process in the asynchronous algorithm, the trajectories are different from each other, resulting in the spread of the trajectories in the transient time. The i.i.d. switching probability $\Pi_i$ that describes asynchronous computing for each distributed

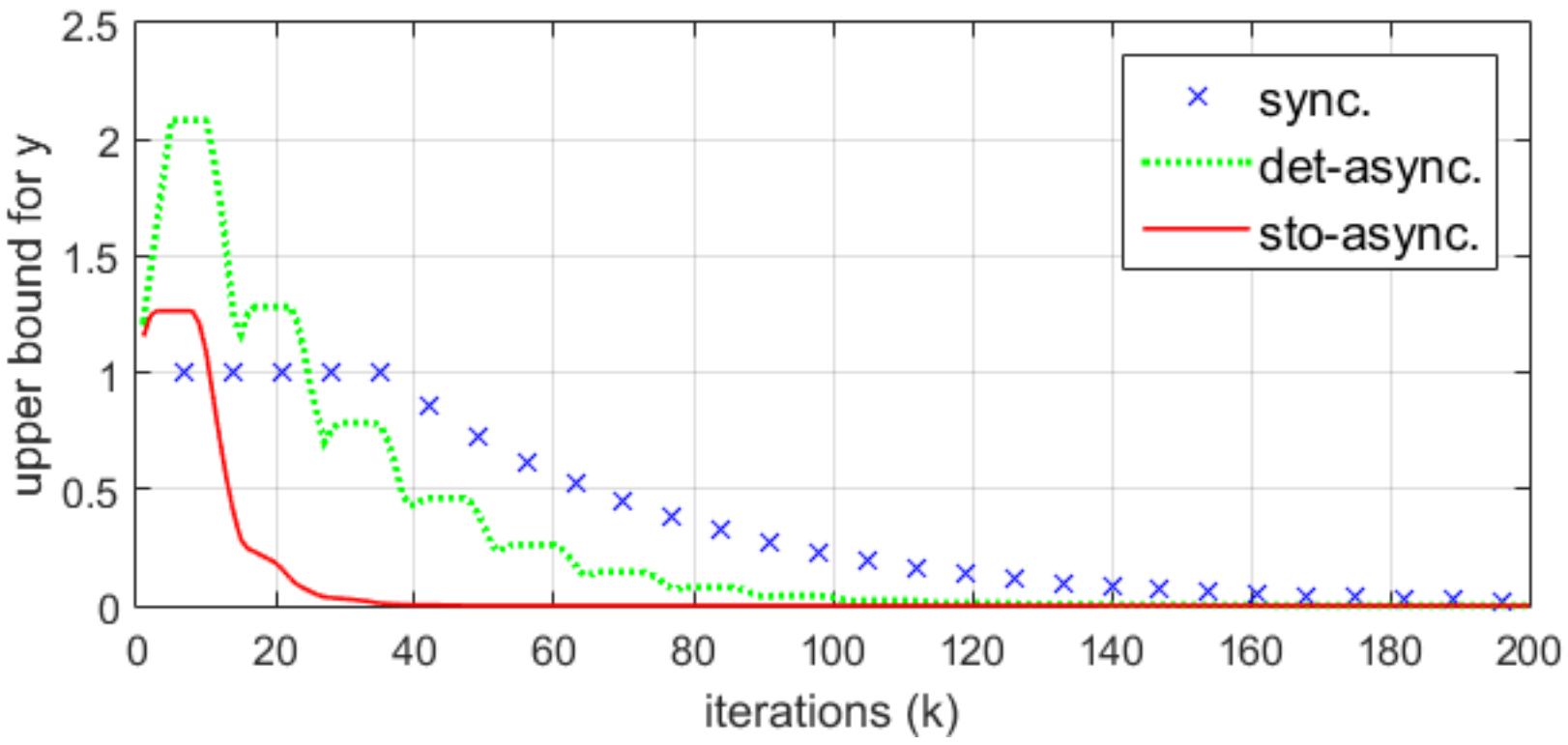

FIG. 4. *The rate of convergence results for distributed quadratic programming with three different schemes: synchronous (cross symbol); deterministic asynchronous (green dotted line); stochastic asynchronous (red solid line) algorithms*

node is given by $(\pi_j)_i = \frac{e^{-3qj}}{\sum_{j=1}^{q} e^{-3qj}}$, $j = 1, 2, \ldots, q$, $\forall i$. Then, by Theorem 4.3 the switching probability for the switched system in (3.2), denoted by $\Pi$, is computed as $\Pi = \big[0,\ 0,\ 0.08,\ 0.8,\ 0.11,\ 0.01,\ 0,\ 0\big]$. For this i.i.d. switching probability, we calculated the spectral radius given in (4.3), which is $\rho\left(\sum_{j=1}^{m} \pi_j (W_j \otimes W_j)\right) = 0.6147 < 1$. Therefore, the convergence of the stochastic asynchronous algorithm is guaranteed in the mean square sense. The result in Fig. 3 also verifies the mean square convergence. The empirical mean and standard deviations are denoted by (red) solid line with cross mark and vertical bars, respectively. As the iteration step increases, the error of the mean square converges to zero (zero standard deviation).

Next, we predict the rate of convergence for three different schemes: i) synchronous case; ii) deterministic asynchronous case; iii) stochastic asynchronous case, in order to compare the performance. By employing the proposed results in section 5, we plotted the rate of convergence in Fig. 4. According to this result for the upper bound of the rate of convergence, the stochastic asynchronous algorithm is advantageous to speedup the total computation time in finding the optimal solution. This stochastic asynchronous scheme is up to 5 times faster than the synchronous algorithm and 2.5 times faster than the deterministic asynchronous algorithm, respectively.

In Fig. 5, we plotted actual computation time to find the optimal solution for three different schemes. For comparison purpose, the computation time for the sequential case is also given as a reference. The termination for the iteration is given by the residual tolerance $|y^k - y^{k-1}| \leq 10^{-5}$. As shown in Fig. 5, the proposed stochastic algorithm achieves the fastest convergence to solve the distributed QP problem. This result coincides with the result on the rate of convergence, which provides information regarding which schemes are the best to solve the given QP problem even before solving the optimization problem.

For three different schemes, Table 1 and Fig. 6 present the computation time and speedup, respectively as we increase the number of threads in the test bed. Also, we plotted speedup of three different schemes based on Table 1, by increasing the total number of threads. As the number of threads increases, the performance degradation occurred in the synchronous case, whereas the deterministic and stochastic asynchronous algorithms resulted in continuous speedup. When the number of threads is 8, the stochastic asynchronous algorithm led to 5.81 times speedup compared to the sequential computing, which is also 2.26 times faster than synchronous algorithms.

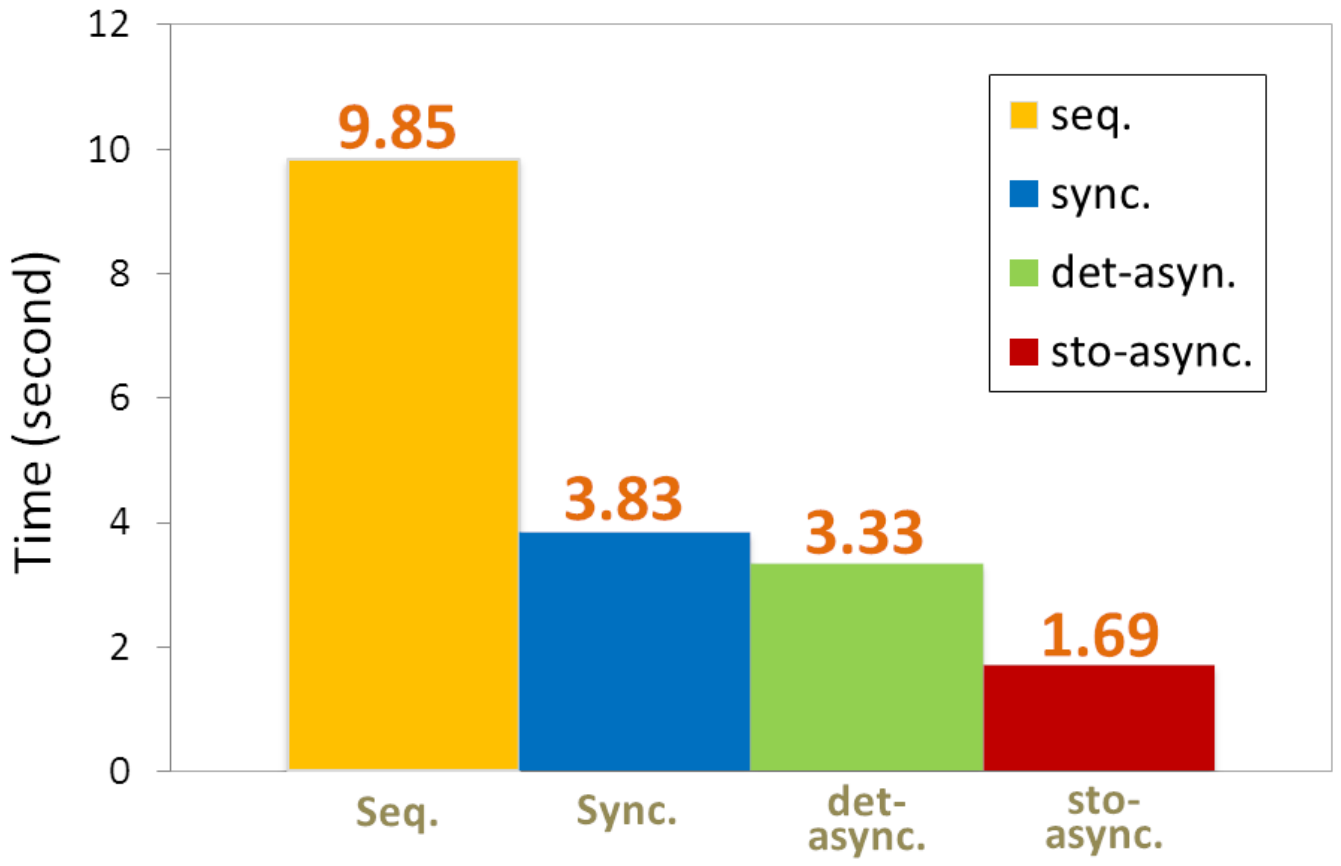

FIG. 5. *The convergence time comparison between sequential computing and three different schemes when the number of threads is given by 8 (maximum possible parallelization for the test bed.*

TABLE 1
*Comparison of total computation time for the dual variable being convergent to the optimal value.*

| No. of Threads | Synchronous | | Det-Asynchronous | | Sto-Asynchronous | |
|---|---|---|---|---|---|---|
| | Time | Speedup | Time | Speedup | Time | Speedup |
| #2 | 5.2012s | 1.89 | 7.8774s | 1.25 | 4.4422s | 2.22 |
| #3 | 4.0189s | 2.45 | 5.8558s | 1.68 | 3.1259s | 3.15 |
| #4 | 3.3848s | 2.91 | 4.8792s | 2.02 | 2.6342s | 3.74 |
| #5 | 3.3511s | 2.94 | 4.3913s | 2.24 | 2.3071s | 4.27 |
| #6 | 3.3547s | 2.94 | 3.8129s | 2.58 | 2.0249s | 4.86 |
| #7 | 3.5891s | 2.74 | 3.4590s | 2.85 | 1.8351s | 5.37 |
| #8 | 3.8340s | 2.57 | 3.3260s | 2.96 | 1.6933s | 5.81 |

As described in Remark 3.1, the computational complexity was the major concern when adopting the switched system framework for analysis of the stochastic asynchronous algorithm. To circumvent this complexity issue, we applied Algorithm 1. Thus, the number of switching modes has been drastically reduced from $q^N = 8^{20000}$ to $q = 8$, owing to Algorithm 1. Consequently, the analysis of stochastic asynchronous computing algorithm was carried out in a computationally efficient manner.

**7. Conclusion.** In this paper, we studied the convergence of asynchronous distributed QP problems via dual decomposition technique. To analyze the behavior of asynchrony in distributed and parallel computing, the switched system framework was introduced. Since the switching mode number becomes large for massively asynchronous computing algorithm, we developed a new algorithm, which drastically decreases mode numbers. By implementing the proposed method, the convergence condition in the mean square sense can be checked without any computational complexity issues. Also, we derived the rate of convergence for three different schemes (e.g., synchronous, deterministic asynchronous, and stochastic asynchronous algorithms), which analytically shows how fast dual variable converges to the optimal solution. The numerical example with an implementation of asynchronous distributed QP using OpenMP supports the validity of the proposed methods.

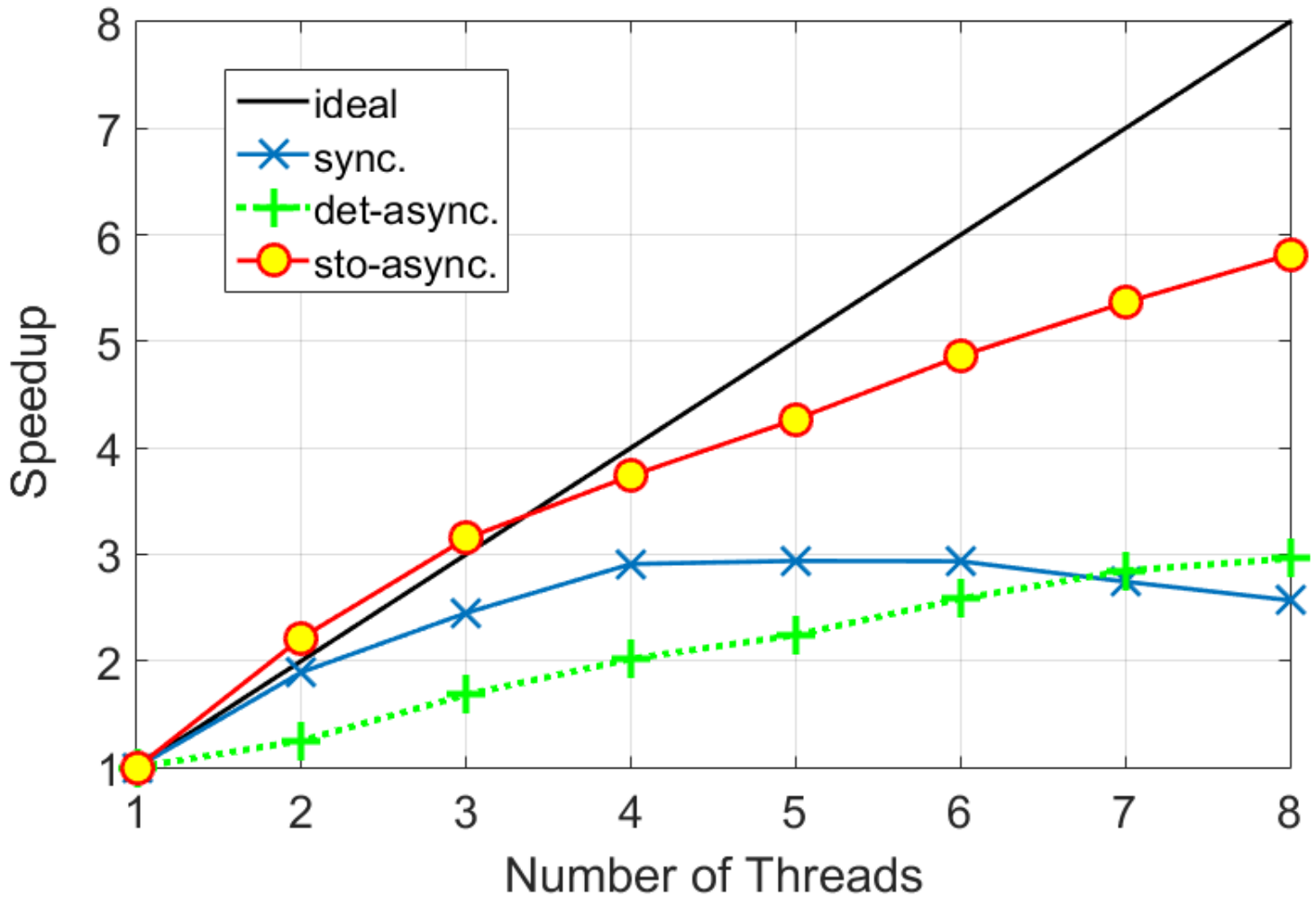

FIG. 6. *The speedup vs. numbers of threads*

**Acknowledgement.** This research material is based upon work supported by the National Science Foundation under grant number #1349100.